\documentclass[10pt]{amsart}
\usepackage{latexsym,amssymb,graphicx,url}

\newcommand{\Inv}{\operatorname{Inv}}
\newcommand{\SYT}{\operatorname{SYT}}
\theoremstyle{plain}
  \newtheorem{theorem}{Theorem}[section]
  \newtheorem{proposition}[theorem]{Proposition}
  \newtheorem{lemma}[theorem]{Lemma}

  \newtheorem*{reform-theorem}{Theorem \ref{han-original}$^\prime$}
\theoremstyle{definition}

  \numberwithin{equation}{section}

\begin{document}

\title{Proof of Han's Hook Expansion Conjecture}

\author{Kevin Carde}
\address{Department of Mathematics\\Harvard University\\Cambridge, MA 02138}
\email{kcarde@fas.harvard.edu}

\author{Joe Loubert}
\address{School of Mathematics\\University of Minnesota\\Minneapolis, MN 55455}
\email{loub0007@umn.edu}

\author{Aaron Potechin}
\address{Department of Mathematics\\Princeton University\\Princeton, NJ 08544}
\email{potechin@princeton.edu}

\author{Adrian Sanborn}
\address{Department of Mathematics\\Harvard University\\Cambridge, MA 02138}
\email{asanborn@fas.harvard.edu}

\date{\today}

\thanks{Special thanks to Victor Reiner and Dennis Stanton for their guidance. This paper was researched and written during their REU program at the University of Minnesota, funded by NSF grants DMS - 0601010 and DMS - 0503660. Also thanks to Andy Manion for his input and friendship throughout the REU}

\begin{abstract}
We prove a conjecture by Guo-Niu Han which interpolates between two known hook expansion formulas. \end{abstract}

\maketitle


\section{Introduction}
\label{intro}
This paper proves a recent conjecture of Guo-Niu Han (\cite[Conjecture 1.4]{Han} and \cite[Conjecture 2.1]{Han2}) giving a hook expansion formula related to partitions, permutations, and involutions.  We refer the reader to \cite[Chapter 7]{Stanley} for the notation used throughout. Han's conjecture interpolates between two hook expansion formulas that follow from classical results on $f^\lambda$, the number of standard Young tableaux (SYT) of shape $\lambda$

The Hook Formula of Frame, Robinson, and Thrall \cite{Frame} states
\begin{equation}
\label{hook-formula}
f^\lambda = \frac{n!}{\displaystyle\prod_{x\in\lambda} h(x)}
\end{equation}
where $h(x)$ denotes the hook length at the cell $x$ in $\lambda$. Elementary representation theory or the Robinson-Schensted-Knuth Algorithm (RSK) (see e.g. \cite[\S3.3]{Sagan}) shows that
\begin{equation*}
\label{rsk}n! = \sum_{\lambda\vdash n} (f^\lambda)^2.
\end{equation*}
Combining these two identities gives us the first hook expansion formula
\begin{equation}\label{hook-exp-1}e^t = \sum_{n=0}^\infty t^n \sum_{\lambda\vdash n} \prod_{x\in\lambda}\frac{1}{h(x)^2}.
\end{equation}

Let $\Inv(n) = \{\pi\in S_n\mid \pi = \pi^{-1}\}$ denote the set of involutions in $S_n$, with the convention that $S_0=\Inv(0)=\{1\}$, where $1$ denotes the identity permutation of the empty set. Elementary representation theory or RSK also shows that
\begin{equation*}
\label{rsk-inv}
|\Inv(n)| = \sum_{\lambda\vdash n} f^\lambda.
\end{equation*}
This gives a second hook expansion formula
\begin{equation}
\begin{split}
\label{hook-exp-2}
e^{t+\frac{t^2}{2}} &= \sum_{n=0}^\infty |\Inv(n)| \frac{t^n}{n!} \\
						&= \sum_{n=0}^\infty t^n \sum_{\lambda\vdash n}\prod_{x\in\lambda}\frac{1}{h(x)}.
\end{split}
\end{equation}
Theorem~\ref{han-original} below, conjectured by Han, yields \eqref{hook-exp-1} upon setting $z=0$ and \eqref{hook-exp-2} upon setting $z=1$.
\begin{theorem}
\label{han-original}
\begin{equation}
\label{han-orig-eqn}e^{t+z\frac{t^2}{2}} = \sum_{n=0}^\infty t^n \sum_{\lambda\vdash n}\prod_{x\in\lambda} \rho(h(x),z)
\end{equation}
where
\begin{equation}
\label{han-weight-function}
\rho(n,z) = \frac{\displaystyle\sum_{k \geq 0} \binom{n}{2k} z^k}{\displaystyle n \sum_{k \geq 0} \binom{n}{2k+1} z^k}.
\end{equation}
\end{theorem}

The left hand side of equation~\ref{han-orig-eqn} has a well-known interpretation as the exponential generating function for involutions counted according to their number of $2$-cycles; see Section~\ref{inv-recur-section} below. We show in Section~\ref{han-equiv-section} below that Theorem \ref{han-original} has the following equivalent reformulation.
\begin{reform-theorem}[Reformulation of Theorem \ref{han-original}]
For all $n\geq 0$,
\begin{equation}
\label{han-new-eqn}
\sum_{\pi\in\Inv(n)} \left(\frac{1+q}{1-q}\right)^{\alpha_1(\pi)} = \sum_{\lambda\vdash n} f^\lambda\prod_{x\in\lambda} \frac{1+q^{h(x)}}{1-q^{h(x)}}
\end{equation}
where $\alpha_1(\pi)$ is the number of fixed points of the permutation $\pi$.
\end{reform-theorem}

\noindent Accordingly, define
\begin{align*}
w(h) &= \frac{1+q^h}{1-q^h}\\
w(\lambda) &= \prod_{x\in\lambda} w(h(x)) = \prod_{x\in\lambda}\frac{1+q^{h(x)}}{1-q^{h(x)}}.
\end{align*}

\noindent In Section~\ref{ext-ret-section}, Theorem~\ref{han-original}$^\prime$ is deduced from the following result, proven in Sections~\ref{ext-ret-proof-section}~and~\ref{multivariate-section}.
\begin{lemma}
\label{ext-ret}
Fix $\lambda\vdash n$. Then
\begin{equation*}
\sum_{\lambda^+\gtrdot\lambda} w(\lambda^+) = w(1) w(\lambda) + \sum_{\lambda^-\lessdot\lambda} w(\lambda^-)
\end{equation*}
where $\lambda^+\gtrdot\lambda$ (resp. $\lambda^-\lessdot\lambda$) indicates that $\lambda^+$ (resp. $\lambda^-$) is obtained by adding (resp. removing) a square to $\lambda$.
\end{lemma}

\section{Generating function and recursion for involutions}
\label{inv-recur-section}
Standard exponential generating function techniques (see e.g. \cite[Eqn. (5.30)]{Stanley}) show the following result due to Touchard, and its consequence for involutions.
\begin{proposition}
\label{Touchard-prop}
If $\alpha_i(\pi)$ denotes the number of $i$-cycles in $\pi$, then
\begin{equation}
\label{Touchard-equation}
\sum_{n=0}^\infty \frac{t^n}{n!}\left( \sum_{\pi\in S_n} u_1^{\alpha_1(\pi)} u_2^{\alpha_2(\pi)} u_3^{\alpha_3(\pi)} \cdots \right)  = e^{u_1\frac{t^1}{1}+u_2 \frac{t^2}{2}+u_3 \frac{t^3}{3}+ \cdots}.
\end{equation}
\end{proposition}

\noindent Setting $u_i = 0$ for $i \geq 3$ yields
\begin{equation}
\label{hanLHS}
e^{u_1t+u_2 \frac{t^2}{2}} = \sum_{n=0}^\infty \frac{t^n}{n!}\sum_{\pi\in\Inv(n)} u_1^{\alpha_1(\pi)} u_2^{\alpha_2(\pi)}.	
\end{equation}
Direct combinatorial reasoning, or differentiation of \eqref{hanLHS} with respect to $t$ gives a well-known recursion for
$$
g_n:=\sum_{\pi\in\Inv(n)} u_1^{\alpha_1(\pi)} u_2^{\alpha_2(\pi)},
$$
namely
\begin{equation}
\label{general-involution-recursion}
g_{n+1} = u_1 g_n + n u_2 g_{n-1}.
\end{equation}

\section{Equivalence of Theorem \ref{han-original} and Theorem \ref{han-original}$^\prime$}
\label{han-equiv-section}

Using the Binomial Theorem, we can rewrite the weight function \eqref{han-weight-function} as
\begin{equation*}
\rho(n,z) = \frac{(1+\sqrt z)^n + (1-\sqrt z)^n}{(1+\sqrt z)^n - (1-\sqrt z)^n}\cdot \frac{\sqrt z}{n} = \frac{1+\left(\frac{1-\sqrt z}{1+\sqrt z}\right)^n}{1-\left(\frac{1-\sqrt z}{1+\sqrt z}\right)^n} \cdot \frac{\sqrt z}{n}.
\end{equation*}

\noindent Starting with Theorem \ref{han-original}, substituting
\begin{align*}
q &= \frac{1-\sqrt z}{1+\sqrt z} \\
T &= t\sqrt z
\end{align*}
and using the Hook Formula \eqref{hook-formula} gives
\begin{equation}
\begin{split}
\label{sub}
e^{\frac{1+q}{1-q}T + \frac{T^2}{2}} &= \sum_{n=0}^\infty \frac{T^n}{z^{n/2}} \sum_{\lambda\vdash n}\prod_{x\in\lambda}\left( \frac{1+q^{h(x)}}{1-q^{h(x)}}\frac{\sqrt z}{h(x)}\right) \\
			&= \sum_{n=0}^\infty \frac{T^n}{n!} \sum_{\lambda\vdash n} f^\lambda\prod_{x\in\lambda}\frac{1+q^{h(x)}}{1-q^{h(x)}} .
\end{split}
\end{equation}
On the other hand, setting $u_1 = (1+q)/(1-q)$ and $u_2 = 1$ in equation \eqref{hanLHS}, we get
\begin{align*}
\sum_{n=0}^\infty \frac{T^n}{n!} \sum_{\pi\in\Inv(n)} \left(\frac{1+q}{1-q}\right)^{\alpha_1(\pi)}  &= e^{\frac{1+q}{1-q}T + \frac{T^2}{2}}\\
			&= \sum_{n=0}^\infty \frac{T^n}{n!} \sum_{\lambda\vdash n} f^\lambda\prod_{x\in\lambda}\frac{1+q^{h(x)}}{1-q^{h(x)}} .
\end{align*}
Equating coefficients of $T^n/n!$ gives Theorem \ref{han-original}$^\prime$.

\section{Lemma~\ref{ext-ret} Implies Theorem \ref{han-original}$^\prime$}
\label{ext-ret-section}

Define
\begin{align*}
\phi_n &= \sum_{\lambda\vdash n} f^\lambda \cdot w(\lambda) \\
\psi_n &= \sum_{\pi\in\Inv(n)} w(1)^{\alpha_1(\pi)}.
\end{align*}
Then Theorem \ref{han-original}$^\prime$ asserts that $\phi_n = \psi_n$ for all $n\geq 0$. By equation~\eqref{general-involution-recursion}, $\psi_n$ satisfies the recursion
$$
\psi_{n+1} = w(1)\psi_n + n\cdot \psi_{n-1}.
$$
Since $\psi_0 = 1 = \phi_0$ and $\psi_1 = w(1) = \phi_1$, it suffices to show that $\phi_n$ satisfies the same recursion, namely
$$
\phi_{n+1} = w(1)\phi_n + n\cdot \phi_{n-1}.
$$

Let $\SYT(n)$ denote the set of all standard Young tableaux of size $n$, and for a tableau $P$, let $\lambda(P)$ be the partition $\lambda$ giving its shape. Notice that one can alternatively express
$$
\phi_n  = \sum_{P\in\SYT(n)} w(\lambda(P)).
$$
Suppose Lemma \ref{ext-ret} is true; i.e., for all $\lambda\vdash n$,
\begin{equation}
\label{ext-ret-eqn}
\sum_{\lambda^+\gtrdot\lambda} w(\lambda^+) = w(1) w(\lambda) + \sum_{\lambda^-\lessdot\lambda} w(\lambda^-)
\end{equation}
where $\mu \gtrdot\lambda$ indicates that $\mu$ is a partition such that $\mu > \lambda$ in the inclusion ordering and $|\mu| = |\lambda| + 1$. Summing \eqref{ext-ret-eqn} over all SYT $P$ of shape $\lambda$ for shapes $\lambda\vdash n$, one obtains
\begin{equation}
\label{Im-uncreative}
\sum_{P\in\SYT(n)}\sum_{\lambda^+\gtrdot\lambda(P)} w(\lambda^+) = w(1)\sum_{P\in\SYT(n)}w(\lambda(P)) + \sum_{P\in\SYT(n)}\sum_{\lambda^-\lessdot\lambda(P)} w(\lambda^-).
\end{equation}
In the sum on the left hand side, we can lift a SYT $P$ of $\lambda$ to a SYT $P^+$ of $\lambda^+$ by labeling the new square in $\lambda^+$ with the number $n+1$. Indeed, every such $P^+$ is clearly obtained exactly once in this way. Thus, \eqref{Im-uncreative} is equivalent to
\begin{equation}
\label{uncreative2}
\sum_{P^+\in\SYT(n+1)} w(\lambda(P^+)) = w(1)\sum_{P\in\SYT(n)}w(\lambda(P)) + \sum_{P\in\SYT(n)}\sum_{x} w(\lambda(P)-x) 
\end{equation}
where the last sum is over corner cells $x\in\lambda(P)$.

We wish to simplify the second term on the right hand side of equation \eqref{uncreative2}. Note that reverse row-insertion on the tableau $P$ starting in the corner cell $x$ produces a tableau $P'$ together with a row-ejected letter $i$. Decrementing by one all entries in $P'$ which are greater than $i$ yields a tableau $P^-$ in $\SYT(n-1)$.This establishes a bijection
$$
\left\{(P,x)\bigg|\begin{array}{c} P\in\SYT(n)\\x\text{ a corner of }P\end{array}\right\}
\longleftrightarrow
\left\{(P^-,i)\bigg|\begin{array}{c}P^-\in\SYT(n-1)\\i\in\{1,\cdots, n\}\end{array}\right\}
$$
which yields the identity
\begin{align*}
\sum_{P\in\SYT(n)}\sum_{x} w(\lambda(P)-x) &= \sum_{P^-\in\SYT(n-1)} \sum_{i=1}^n w(\lambda(P^-)).\\
																					 &= n\sum_{P^-\in\SYT(n-1)} w(\lambda(P^-))
\end{align*}
Finally, substituting this into \eqref{uncreative2} gives
\begin{equation*}
\sum_{P^+\in\SYT(n+1)} w(\lambda(P^+)) = w(1)\sum_{P\in\SYT(n)}w(\lambda(P)) + n\sum_{P^-\in\SYT(n-1)} w(\lambda(P^-)) 
\end{equation*}
which gives the desired recursion
\begin{equation*}
\phi_{n+1} = w(1)\phi_n + n\cdot\phi_{n-1}.
\end{equation*}

\section{Proof of Lemma~\ref{ext-ret}}
\label{ext-ret-proof-section}

\begin{figure}[h]
\includegraphics[width=1\textwidth]{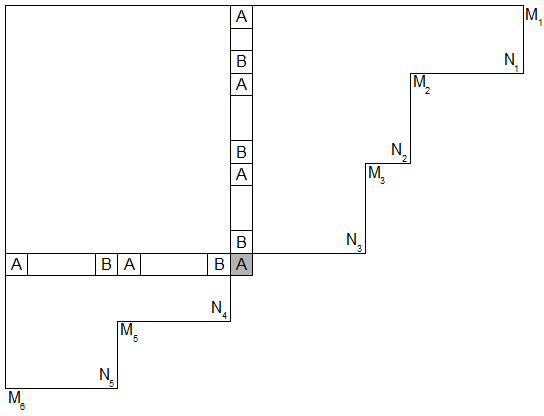}
\caption{\label{ext-ret-example} Adding an outer corner at $M_4$. Hooks at squares labeled `A' (resp. `B') remain uncanceled in the numerator (resp. denominator).}
\end{figure}

For a partition $\lambda$, label the outer corners of $\lambda$ as $M_{1},...,M_{d}$ with coordinates $(a_1,b_1),...,(a_{d},b_{d})$ and label the inner corners (i.e., $1$-hooks) as $N_{1},...,N_{d-1}$ with coordinates $(\alpha_1,\beta_1),...,(\alpha_{d-1},\beta_{d-1})$; see Figure~\ref{ext-ret-example} for an example. Define the {\it content} of the square $(i,j)$ to be $c(i,j)=j-i$. To prove Lemma~\ref{ext-ret} for $\lambda$, we will reduce it to an equation relating the contents of the inner and outer corners\footnote{Ideas used in this proof go back to a proof of the Hook Formula by Vershik \cite{Vershik}, also investigated by Kirillov \cite{Kirillov} and Kerov \cite{Kerov1, Kerov2}.} of $\lambda$. Then we will prove that this equation is in fact true with the contents replaced by arbitrary variables.

If $\lambda^{+}$ is obtained from $\lambda$ by adding an outer corner $M_{k}$, then we can find an explicit formula for $w(\lambda^{+})$ in terms of $w(\lambda)$ and the contents of the outer and inner corners. The terms of $w(\lambda^{+})$ mostly agree with the terms of $w(\lambda)$ because the hook length of a square will only change if it is in the same row or column as $M_{k}$. Due to the changes in hook lengths at a square, we introduce the notation $h_\lambda(i,j)$ for the hook length at the square in row $i$ and column $j$ of the shape $\lambda$. If the hook length at a square changes in passing from $\lambda$ to $\lambda^+$, it must increase by one. Finally, we introduce the $1$-hook $M_{k}$ as an extra factor in $w(\lambda^+)$:
\begin{equation*}
\frac{w(\lambda^{+})}{w(\lambda) w(1)}    = \prod_{j=1}^{a_{k}-1}\frac{w(h_{\lambda^+}(j,b_{k}))}{w(h_{\lambda}(j,b_{k}))}      \prod_{j=1}^{b_{k}-1}\frac{w(h_{\lambda^+}(a_{k},j))}{w(h_{\lambda}(a_{k},j))}.
\end{equation*}

Within these products, more terms cancel. If there is no inner corner in row $j$, then $h_{\lambda}(j,b_{k})=h_{\lambda}(j+1,b_{k})+1$. Also, row $j$ has an inner corner if and only if row $j+1$ has an outer corner. Hence, if row $j+1$ has no inner corner, then the term $w(h_{\lambda^+}(j,b_k))$ in the numerator cancels with the term $w(h_{\lambda}(j+1,b_k))$ in the denominator, and only the terms in rows or columns with inner or outer corners remain; see Figure \ref{ext-ret-example} for an example.

This allows us to write
\begin{equation*}
\prod_{j=1}^{a_{k}-1}\frac{w(h_{\lambda^+}(j,b_{k}))}{w(h_{\lambda}(j,b_{k}))} = \frac{w(h_{\lambda^+}(a_{1},b_{k}))}{w(h_{\lambda}(\alpha_{1},b_{k}))} \cdot \frac{w(h_{\lambda^+}(a_{2},b_{k}))}{w(h_{\lambda}(\alpha_{2},b_{k}))} \cdot \cdot \cdot \frac{w(h_{\lambda^+}(a_{k-1},b_{k}))}{w(h_{\lambda}(\alpha_{k-1},b_{k}))}.
\end{equation*}
Similarly, swapping rows and columns above,
\begin{equation*}
\prod_{j=1}^{b_{k}-1}\frac{w(h_{\lambda^+}(a_{k},j))}{w(h_{\lambda}(a_{k},j))} = \frac{w(h_{\lambda^+}(a_{k},b_{d}))}{w(h_{\lambda}(a_{k},\beta_{d-1}))} \cdot \frac{w(h_{\lambda^+}(a_{k},b_{d-1}))}{w(h_{\lambda}(a_{k},\beta_{d-2}))} \cdot \cdot \cdot \frac{w(h_{\lambda^+}(a_{k},b_{k+1}))}{w(h_{\lambda}(a_{k},\beta_{k}))}.
\end{equation*}

The length of the hook of $(a_{i},b_{k})$ is just $c(a_{i},b_{i}-1)-c(a_{k}-1,b_{k})+1 = c(M_i)-c(M_k)-1$, so setting $x_{i}=c(M_{i})$ and $y_{i}=c(N_{i})$,
\begin{equation}
\begin{split}
\label{hc}
h_{\lambda^+}(a_i,b_k) &= x_i-x_k\text{\;\;\;for $i \in \{1,...,k-1\}$} \\
h_{\lambda}(\alpha_i,b_k) &= y_i-x_k\text{\;\;\;for $i \in \{1,...,k-1\}$} \\
h_{\lambda^+}(a_k,b_i) &= x_k-x_i\text{\;\;\;for $i \in \{k+1,...,d\}$} \\
h_{\lambda}(a_k,\beta_i) &= x_k-y_i\text{\;\;\;for $i \in \{k,...,d-1\}$}.
\end{split}
\end{equation}
Then
\begin{equation*}
\frac{w(\lambda^{+})}{w(\lambda) w(1)}  = \frac{\displaystyle \prod_{i=1}^{k-1}w(x_i-x_k)}{\displaystyle \prod_{i=1}^{k-1}w(y_i-x_k)} \,\,    \frac{\displaystyle \prod_{i=k+1}^{d}w(x_k-x_i)}{\displaystyle \prod_{i=k}^{d-1}w(x_k-y_i)}.
\end{equation*}
Summing over $k \in \{1,...,d\}$ we obtain
\begin{equation}
\label{lambda-plus}
\sum_{\lambda^+ \gtrdot \lambda}\frac{w(\lambda^{+})}{w(\lambda) w(1)}\; = \; \sum_{k=1}^d \; \frac{\displaystyle \prod_{i=1}^{k-1}w(x_i-x_k)}{\displaystyle \prod_{i=1}^{k-1}w(y_i-x_k)}     \,\, \frac{\displaystyle \prod_{i=k+1}^{d}w(x_k-x_i)}{\displaystyle \prod_{i=k}^{d-1}w(x_k-y_i)}.
\end{equation}

\noindent Now if $k \in \{1,...,d-1\}$, let $\lambda^-$ be the partition obtained by removing the corner $N_k$ from $\lambda$. A similar formula holds for $w(\lambda^-)$. Again, the only hooks affected by deleting $N_k$ come from squares in the same row or column as $N_k$, giving the equality
\begin{equation*}
\frac{w(\lambda^-) w(1)}{w(\lambda)}  =  \prod_{j=1}^{\alpha_k-1}\frac{w(h_{\lambda^-}(j,\beta_k))}{w(h_{\lambda}(j,\beta_k))}       \prod_{j=1}^{\beta_k-1}\frac{w(h_{\lambda^-}(\alpha_k,j))}{w(h_{\lambda}(\alpha_k,j))}.
\end{equation*}
Again, many of these terms cancel, reducing to
\begin{equation*}
\prod_{j=1}^{\alpha_k-1}\frac{w(h_{\lambda^-}(j,\beta_k))}{w(h_{\lambda}(j,\beta_k))} = \frac{w(h_{\lambda^-}(\alpha_{1},\beta_k)}{w(h_{\lambda}(a_{1},\beta_k))} \cdots  \frac{w(h_{\lambda^-}(\alpha_{k-1},\beta_k)}{w(h_{\lambda}(a_{k-1},\beta_k))} \cdot \frac{w(1)}{w(h_{\lambda}(a_k,\beta_k))}.
\end{equation*}
and
\begin{equation*}
\prod_{j=1}^{\beta_k-1}\frac{w(h_{\lambda^-}(\alpha_k,j))}{w(h_{\lambda}(\alpha_k,j))} = \frac{w(h_{\lambda^-}(\alpha_{k},\beta_{d-1}))}{w(h_{\lambda}(\alpha_{k},b_d))} \cdots  \frac{w(h_{\lambda^-}(\alpha_{k},\beta_{k+1})}{w(h_{\lambda}(\alpha_{k},b_{k+2}))} \cdot \frac{w(1)}{w(h_{\lambda}(\alpha_k,b_{k+1}))}.
\end{equation*}
Analogous to equations \eqref{hc}, we have
\begin{align*}
h_{\lambda^-}(\alpha_i,\beta_k) &= y_i-y_k\text{\;\;\;for $i \in \{1,...,k-1\}$} \\
h_{\lambda}(a_i,\beta_k) &= x_i-y_k\text{\;\;\;for $i \in \{1,...,k\}$} \\
h_{\lambda^-}(\alpha_k,\beta_i) &= y_k-y_i\text{\;\;\;for $i \in \{k+1,...,d-1\}$} \\
h_{\lambda}(\alpha_k,b_i) &= y_k-x_i\text{\;\;\;for $i \in \{k+1,...,d\}$}.
\end{align*}
These allow us to write
\begin{equation*}
\frac{w(\lambda^-)}{w(1)w(\lambda)}  = \frac{\displaystyle \prod_{i=1}^{k-1}w(y_i-y_k)}{\displaystyle \prod_{i=1}^{k}w(x_i-y_k)} \,\,     \frac{\displaystyle \prod_{i=k+1}^{d-1}w(y_k-y_i)}{\displaystyle \prod_{i=k+1}^{d}w(y_k-x_i)}.\end{equation*}
Summing this over $k \in \{1,...,d-1\}$ we have
\begin{equation}
\label{lambda-minus}
\sum_{\lambda^- \lessdot \lambda} \frac{w(\lambda^-)}{w(1)w(\lambda)}  \; = \; \sum_{k=1}^{d-1} \; \frac{\displaystyle \prod_{i=1}^{k-1}w(y_i-y_k)}{\displaystyle \prod_{i=1}^{k}w(x_i-y_k)} \,\,   \frac{\displaystyle \prod_{i=k+1}^{d-1}w(y_k-y_i)}{\displaystyle \prod_{i=k+1}^{d}w(y_k-x_i)}.
\end{equation}
Plugging \eqref{lambda-plus} and \eqref{lambda-minus} into Lemma \ref{ext-ret} and employing the fact that 
$$
w(-x)=-w(x),
$$ 
we are reduced to proving
\begin{proposition}
\label{cont-prop}
For distinct values $x_1,x_2,...,x_d,y_1,y_2,...,y_{d-1}$, one has
\begin{equation}
\label{cont-eq}
\sum_{k=1}^d\frac{\displaystyle \prod_{i=1,i\ne k}^dw(x_k-x_i)}{\displaystyle \prod_{i=1}^{d-1}w(x_k-y_i)}+\sum_{k=1}^{d-1}\frac{\displaystyle \prod_{i=1,i\ne k}^{d-1}w(y_k-y_i)}{\displaystyle \prod_{i=1}^dw(y_k-x_i)}=1.
\end{equation}
\end{proposition}

\noindent Note that we are applying this proposition in the special case where $x_i = c(M_i)$ and $y_i = c(N_i)$, so that $x_1,x_2,...,x_d,y_1,y_2,...,y_{d-1}$ are indeed distinct. Proposition~\ref{cont-prop} is a special case of the following proposition, given two proofs in Section~\ref{multivariate-section}.
\begin{proposition}
\label{cont-lemma}
Within the field of rational functions ${\mathbf{Q}}(a_1,a_2,\ldots,a_n)$, one has
\begin{equation}
\label{cont-lemma-eq}
\sum_{k=1}^n\prod_{i=1,i\ne k}^n\frac{a_k+a_i}{a_k-a_i}=\left\{\begin{array}{l l}0 & \quad\mbox{if $n$ is even}\\1 & \quad\mbox{if $n$ is odd}\end{array}\right..
\end{equation}
\end{proposition}

We now explain how Proposition~\ref{cont-lemma} implies Proposition~\ref{cont-prop}.
For $d=1$ the statement is trivial, so assume $d\geq 2$. Rewrite Proposition \ref{cont-lemma} with $n=2d-1$ as:
\begin{equation*}
\sum_{k=1}^d\frac{\displaystyle\prod_{i=1,i\ne k}^d\frac{a_k+a_i}{a_k-a_i}}{\displaystyle\prod_{i=d+1}^{2d-1}\frac{a_k-a_i}{a_k+a_i}}+\sum_{k=d+1}^{2d-1}\frac{\displaystyle\prod_{i=d+1,i\ne k}^{2d-1}\frac{a_k+a_i}{a_k-a_i}}{\displaystyle\prod_{i=1}^d\frac{a_k-a_i}{a_k+a_i}}=1.
\end{equation*}
Multiply each factor by $a_k^{-1}/a_k^{-1}$:
\begin{equation*}
\sum_{k=1}^d\frac{\displaystyle\prod_{i=1,i\ne k}^d\frac{1+a_ia_k^{-1}}{1-a_ia_k^{-1}}}{\displaystyle\prod_{i=d+1}^{2d-1}\frac{1-a_ia_k^{-1}}{1+a_ia_k^{-1}}}+\sum_{k=d+1}^{2d-1}\frac{\displaystyle\prod_{i=d+1,i\ne k}^{2d-1}\frac{1+a_ia_k^{-1}}{1-a_ia_k^{-1}}}{\displaystyle\prod_{i=1}^d\frac{1-a_ia_k^{-1}}{1+a_ia_k^{-1}}}=1.
\end{equation*}
Now we set $a_i = q^{-x_i}$ for $1\le i\le d$ and $a_i = -q^{-y_{i-d}}$ for $d+1\le i\le 2d-1$:
\begin{equation*}\sum_{k=1}^d\frac{\displaystyle\prod_{i=1,i\ne k}^d\frac{1+q^{x_k-x_i}}{1-q^{x_k-x_i}}}               {\displaystyle\prod_{i=1}^{d-1}\frac{1+q^{x_k-y_i}}{1-q^{x_k-y_i}}}+\sum_{k=1}^{d-1}\frac{\displaystyle\prod_{i=1,i\ne k}^{d-1}\frac{1+q^{y_k-y_i}}{1-q^{y_k-y_i}}}{\displaystyle\prod_{i=1}^d\frac{1+q^{y_k-x_i}}{1-q^{y_k-x_i}}}=1.
\end{equation*}
This is precisely equation~\eqref{cont-eq} upon plugging in $\displaystyle w(h)=\frac{1+q^h}{1-q^h}$.

Assuming for the moment Proposition~\ref{cont-lemma}, this establishes Lemma~\ref{ext-ret}, which proves Theorem \ref{han-original}$^\prime$ and Theorem \ref{han-original}.

\newpage

\section{Two proofs of Proposition \ref{cont-lemma}}
\label{multivariate-section}

\noindent
\emph{Proof 1: }
Set
\begin{equation*}
b_k:=\prod_{i=1,i\ne k}^n\frac{a_k+a_i}{a_k-a_i}.
\end{equation*}
We wish to show that the sum of the $b_k$ is $0$ or $1$ depending on the parity of $n$.

Consider the partial fraction decomposition
\begin{equation}
\label{t-frac}\prod_{i=1}^n\frac{t+a_i}{t-a_i}=c_0+\sum_{k=1}^n\frac{c_k}{t-a_k}.
\end{equation}
Taking the limit $t\to\infty$ on both sides yields $c_0=1$. Multiplying both sides by $t-a_k$ and setting $t=a_k$ gives
\begin{equation*}
2a_k\prod_{i=1,i\ne k}^n\frac{a_k+a_i}{a_k-a_i}=c_k.
\end{equation*}
So $c_k=2a_kb_k$. Setting $t=0$ in \eqref{t-frac} gives
\begin{equation*}
(-1)^n=c_0-\sum_{k=1}^n\frac{c_k}{a_k}.
\end{equation*}
Plugging in $c_0 = 1$ and $c_k = 2a_k b_k$ yields
\begin{equation*}
1-(-1)^n=2\sum_{k=1}^nb_k.
\end{equation*}
The left hand side is $0$ if $n$ is even and $2$ if $n$ is odd, so dividing through by $2$ yields the desired result.\qed
\vskip.2in
\noindent
\emph{Proof 2: }
Multiply through by the denominator in \eqref{cont-lemma-eq}, so that the equation to be proved is
\begin{equation}
\label{cont-lemma-joe}
\sum_{k=1}^n(-1)^{k-1}\prod_{\substack{i=1\\i \ne k}}^n(a_k+a_i)\prod_{\substack{i<j \\ i \ne k \\ j \ne k}}(a_i-a_j) = \delta_n \cdot \prod_{i<j}(a_i-a_j)
\end{equation}
where $\delta_n$ is 0 if $n$ is even and $1$ if $n$ is odd.

We wish to show that the polynomial on the left hand side of equation \eqref{cont-lemma-joe} is an alternating function of the variables $a_1,\cdots,a_n$.  Consider the effect of exchanging the variables $a_r, a_{r+1}$.
For $k \ne r,r+1$, the only change in the summand is that the $a_r-a_{r+1}$ in the second product is replaced with $a_{r+1}-a_r$, changing the sign. For $k=r,r+1$, the summand itself stays the same, but the $(-1)^{k-1}$ factor is off by one on each summand, again changing the sign, as desired. Since the left hand side of equation \eqref{cont-lemma-joe} is alternating, the Vandermonde product in the right hand side divides the left hand side.  The left hand side has degree at most $\binom{n}{2}$.
Therefore, equation \eqref{cont-lemma-joe} holds for some constant $\delta_n$, which is determined from consideration of the coefficient of $a_1^{n-1}a_2^{n-2} \cdots a_{n-1}$  on each side.\qed \\

\end{document}